\documentstyle[12pt]{article}
\begin{document}
\newtheorem{lem}{Lemma}
\newtheorem{rem}{Remark}
\newtheorem{question}{Question}
\newtheorem{prop}{Proposition}
\newtheorem{cor}{Corollary}
\newtheorem{thm}{Theorem}
\newtheorem{definition}{Definition}
\newtheorem{exam}{Example}
\title
{A miscellany}
\author{
YanYan Li\thanks{Partially
 supported by
        NSF grant  DMS-0701545.}\\
Department of Mathematics\\
Rutgers University\\
110 Frelinghuysen Road\\
Piscataway, NJ 08854\\
\\
Louis Nirenberg
 \\
Courant Institute\\
251 Mercer Street \\
New York, NY 10012\\
}

%\date{ }
\maketitle

%%%%%  The following two lines should be changed
%\input amstex
\input { amssym.def}
%%%%%%%%%

\centerline{\it Dedicated to the memory of
Vittorio Cafagna}

\bigskip

\bigskip

\setcounter {section} {0}

Vittorio and I, the second named
author, were friends for many years.  We first met at a conference
in Bari in 1975 where
we discussed some work he was doing.  It was only several years
later that we really became
friends {\bf --- } through Henri Berestycki who was
staying with Vittorio 
and his wife in New York.

After that we would meet in Italy
or in New York, and it was always a great pleasure to spend time
together.   We would discuss problems on which he was working ---
they were always interesting, and I was infected by his enthusiasm 
{\bf --- } as well as everything 
under the sun: politics, music, etc.  I well recall one evening
we spent at a performance of tango.

Vittorio was full of life and seemed to enjoy
whatever came his way.  The last time we met was in 
Napoli, several years ago.  He came from Salerno; the evening was 
delightful.

That$'$s how my memory is of him: full of
fun, and interest in everything.

In recent years he turned to use of mathematics in music, with his usual
contagious enthusiasm.

His life was cut short too soon.

\bigskip

In this paper we present several results, unrelated
to each other.  

\section{Extension of a simple inequality}

In \cite{L-N1} 
we gave some extensions of the following well known inequality:

\begin{prop}
Let $u\in C^2(-R, R)$, $u\ge 0$ and
assume 
\begin{equation}
|\ddot u|\le M.
\label{1.1}
\end{equation}
Then
\begin{equation}
|\dot u(0)|\le \sqrt{  2u(0)M}\ \
\qquad\ \ \  \mbox{if}\ M\ge \frac {2u(0)}{R^2},
\label{1.2}
\end{equation}
\begin{equation}
|\dot u(0)|\le \frac{u(0)}R+\frac R2 M\ \
\qquad \mbox{if}\ M<  \frac {2u(0)}{R^2}.
\label{1.3}
\end{equation}
\label{proposition1}
\end{prop}

In \cite{L-N1} this was extended to functions
$u\ge 0$ in higher dimensions, with
(\ref{1.1}) replaced by
\begin{equation}
|\Delta u|\le M.
\label{1.4}
\end{equation}

In \cite{CD-V}, I. Capuzzo Delcetta and A. Vitola 
give far reaching
extensions of the result to fully nonlinear
elliptic operators.

\medskip

Here are some further remarks.

\medskip

\noindent{\bf Proposition 
\ref{proposition1}$'$.}\ {\it  Proposition
\ref{proposition1}    still holds
if (\ref{1.1}) is replaced by
\begin{equation}
\ddot u\le M.
\label{1.5}
\end{equation}
}

\medskip

Indeed, if $M=0$, (\ref{1.3}) is
a standard inequality for concave functions.

The proof of  Proposition
\ref{proposition1}$'$ is just the same as that of   Proposition
\ref{proposition1}  --- see for example 
 \cite{L-N1} and  \cite{NT}.

\begin{rem}
In the result of \cite{L-N1} the condition
(\ref{1.4}) above may not be replaced by
$\Delta u\le M$.
\end{rem}

We now give an extension of  Proposition
\ref{proposition1} to the heat operator.

In the interesting papers \cite{M-K} and 
\cite{M-S} various other extensions of Proposition \ref{proposition1}
have been proved.

We consider $u\ge 0$, $u(t,x)$,
in a parabolic cylinder
$$
B_R=\{(x, t)\ |\ -R^2\le t\le 0, \ x\in \Bbb R^n,\ 
|x|<R\},
$$
satisfying
\begin{equation}
|u_t-\Delta u|\le M.
\label{1.6}
\end{equation}
Here $\Delta$ is the Laplace operator in the $x-$variables.
The origin, $(0,0)$, plays the role of the center of the
cylinder.
\begin{rem}
In case $M=0$, i.e. $u$ satisfies the 
heat equation, but it is not true that
\begin{equation}
|\nabla_x u(0, 0)|\le \frac CR u(0, 0).
\label{1.7}
\end{equation}
\label{remark2}
\end{rem}
Indeed, for $R=1$ and $n=1$, take
$$
u=\frac 1{   (t+1)^{ 1/2}  } e^{ -\frac {  (x+\xi)^2}  { 4(t+1) }  }.
$$
Then
$$
u_x(0, 0)=-\frac \xi 2 u(0, 0).
$$
$\xi$ can be arbitrarily large, (\ref{1.7}) cannot hold.
Instead of (\ref{1.7}) we can bound $|\nabla _x u|$ in a suitable subregion
of $B_R$:  For $r\le R$ consider
$$
A_r=\{ (t, x)\ |\
-\frac 23 r^2 <t <-\frac 13 r^2,\
|x|<\frac r{\sqrt{3} }\}\subset B_{ \sqrt{ \frac 23}R}.
$$

\begin{thm}
Suppose $u\ge 0$ and
\begin{equation}
|u_t-\Delta u|\le M\qquad
\mbox{in}\ B_R.
\label{1.8}
\end{equation}
Then, 
\begin{equation}
|\nabla_x u|+\frac ur \le \overline C
\sqrt{  u(0, 0)M}\quad
\mbox{in}\ A_r\ \mbox{if}\
r=\sqrt{ \frac{u(0, 0)}M}\le R,
\label{1.9}
\end{equation}
\begin{equation}
|\nabla_x u|+\frac uR \le \overline C
\left( \frac{u(0, 0)}R +MR\right)\quad \mbox{in}\ 
A_R\  \mbox{if}\
R< \sqrt{ \frac{u(0, 0)}M}.
\label{1.10}
\end{equation}
Here $\overline C$ depends only on $n$.
\label{theorem1}
\end{thm}

In the elliptic case in \cite{L-N1}, the proof 
was based on standard elliptic estimates and
the Harnack inequality.  Here we make use of the
 corresponding Harnack inequality for nonnegative solutions 
$v$ of the heat equation.  It was proved independently
in 1954 by Hadamard \cite{H}, Pini \cite{P}; see
 Evans \cite{E}.  
We formulate it in a form convenient for our application.

\medskip

\noindent
\underline{Parabolic Harnack Inequality}.\
Assume
$$
v_t-\Delta v=0, \quad v>0\qquad \mbox{in}\ B_\rho.
$$
Then 
$$
\max_{V_\rho}
 v\le C\min_{|x|\le \frac 89\rho}v(x, 0),
$$
where
$$
V_\rho=\{(x, t)\ |\ |x|\le  \frac 89\rho,
\ -\frac 89\rho^2\le t\le -\frac {\rho^2}9\}
$$
and
 $C$ is independent of $\rho$.

In particular, it follows that 
$$
\max_{V_\rho} v\le Cv(0,0).
$$

\medskip

Using the parabolic Harnack inequality and some
standard estimates for the
heat operator, we now present the 

\noindent\underline{Proof of Theorem \ref{theorem1}}.\
For $0<r<R$ let $v$ be the 
solution of 
\begin{eqnarray*}
v_t-\Delta v&=&0 \qquad\mbox{in}\ B_r,
\\
v&=&u
\qquad \mbox{on}\ \partial_p B_r.
\end{eqnarray*}
where $\partial_p B_r$ is the parabolic boundary of
$B_r$:
$$
\partial_p B_r= \partial B_r\setminus\{t=0\}.
$$

We use the standard inequality:
\begin{equation}
v+r|\nabla_x v|\le C_1\max_{ D_r} v\qquad\mbox{in}\ A_r.
\label{1.12}
\end{equation}
Here
$$
D_r=\{(t,x)\ |\ -\frac 56 r^2 <t<
- \frac{r^2}3,\ |x|<\frac r{ \sqrt{2} }\}
$$
and $C_1$ is independent of $r$.

Now applying the parabolic Harnack inequality above, with
$
\rho=r,
$
and  we find, using $D_r\subset V_r$,  that
$$
\max_{D_r}v\le C v(0,0).
$$
Combining this with (\ref{1.12}) we obtain:
\begin{equation}
v+r|\nabla_x v|\le CC_1v(0,0)\qquad\mbox{in}\ A_r.
\label{1.13}
\end{equation}
Next, the function
$$
w=u-v\quad\mbox{in}\ B_r
$$
satisfies
$$
|w_t-\Delta w|\le M,\quad w=0\ \ \mbox{on}\
\partial_p B_r.
$$
Using another standard inequality for the heat operator
we have
\begin{equation}
r|\nabla _xw|+|w|\le C_2Mr^2\quad
\mbox{in}\ B_{3r/4} \ (\mbox{containing}\ A_r).
\label{1.14}
\end{equation}
This, together with (\ref{1.13}), yield: in $A_r$,
\begin{eqnarray*}
u+r|\nabla_x u|&\le & 
v+w+r|\nabla_x w|+r|\nabla_x v|\\
&\le & C_2 Mr^2 +CC_1 v(0,0)\\
&\le & C_2 Mr^2 +CC_1 u(0,0)+CC_1 |w(0,0)|\\
&\le& C_2 Mr^2 +CC_1 u(0,0)+CC_1C_2 Mr^2;
\end{eqnarray*}
the last, by (\ref{1.14}) again.  Thus, in $A_r$ we have
\begin{equation}
\frac ur +|\nabla_x u|\le C_3
\left( \frac{u(0,0)}r +Mr\right).
\label{1.15}
\end{equation}

Next,

\noindent\underline{Case 1}.\  
$\displaystyle{
R\ge \sqrt{   \frac{u(0,0) }M  }.
}$

In this case, with
$$
r=  \sqrt{   \frac{u(0,0) }M  },
$$
(\ref{1.15}) yields (\ref{1.9}).

\medskip

\noindent\underline{Case 2}.\
$\displaystyle{
R< \sqrt{   \frac{u(0,0) }M  }.
}$

In this case, (\ref{1.15}) yields (\ref{1.10})
if $r=R$.

\vskip 5pt
\hfill $\Box$
\vskip 5pt

\begin{rem} Theorem \ref{theorem1} also holds if (\ref{1.6}) is replaced
by 
$$
|u_t-\partial_i(a_{ij}(t,x)u_j)|\le M
$$
if $\{a_{ij}\}$ is uniformly positive definite, and smooth.
\end{rem}

\section{On the Hopf Lemma}

In this section we present an extension of the Hopf Lemma.

A standard form of the lemma for nonlinear second order elliptic
operators:
$$
F(x,u,\nabla u, \nabla^2 u)
$$
is the following.  Here we assume 
$F(x, s, p,A)$ is $C^1$, $s\in \Bbb R$,
$p\in \Bbb R^n$, $A\in {\cal S}^{n\times n}$, the
set of $n\times n$ real symmetric matrices,
and strongly elliptic:
$$
\frac{\partial F}{\partial A_{ij}}
\ \mbox{is positive definite}.
$$

\noindent\underline{Hopf Lemma}:\ Suppose $u, v$ are $C^2$ functions
in a domain $\Omega$ in $\Bbb R^n$ with $C^2$ boundary, both continuous
in $\Omega\cup \{y\}$, $y$ a boundary
point, and that
\begin{eqnarray}
u&>&v\qquad\mbox{in}\ \Omega, \label{2.1}\\
u(y)&=&v(y).
\label{2.2}
\end{eqnarray}
Assume that
\begin{equation}
F(x, u, \nabla u, \nabla^2 u)\le F(x,v,\nabla v,\nabla^2 v)\qquad\mbox{in}\
\Omega.
\label{2.3}
\end{equation}
Then, if $\nu$ is the interior unit normal to
$\partial \Omega$ at $y$,
\begin{equation}
\liminf_{ t\to 0^+}
\frac{  (u-v)(y+t\nu) }t>0.
\label{2.4}
\end{equation}

\medskip

\begin{rem}
The result  
holds if $\partial F/\partial u_{ij}$ is positive semidefinite provided
\begin{equation}
\frac {\partial F}{ \partial u_{ij} }\nu_i\nu_j>0
\label{2.5}
\end{equation}
(say, for all values of the arguments in $F$).
\end{rem}

\begin{rem}
Since $\partial \Omega\in C^2$ is assumed in the result, 
one only needs to prove it for $\Omega$ being an open ball
(working with a ball inside $\Omega$ and having
$y$ on its boundary).
\end{rem}

In an unpublished manuscript, \cite{L-N5},
we extended the Hopf Lemma to domains
with $C^{1, \alpha}$ boundary, $0<\alpha<1$.  The paper
was not published because
we learned that the result had been proved earlier by Kamynin  and
 Khimchenko in
 \cite{KK}.

In the last 25 years or 
so many authors have studied various classes of nonlinear 
degenerate elliptic operators (see for example 
\cite{CLN1}, \cite{CLN2}, \cite{C-N-S1},  
\cite{C-N-S},   \cite{H-L},
 \cite{L1},
\cite{L2}
 where many references may be found).
One way of looking at such 
operators, for a $C^2$ function
$u$ in $\Omega$ in $\Bbb R^n$ is to consider
the symmetric matrix 
$$
\widetilde A_u:= \nabla^2 u+L(\cdot, u, \nabla u),
$$
where
 $L\in C^1(\Omega\times \Bbb R\times \Bbb R^n)$,
 is in ${\cal S}^{n\times n}$,
 the set of  $n\times n$ real  symmetric
matrices.

For every $x$ in $\Omega$, the matrix is required to lie
in a region ${\cal G}$ in 
${\cal S}^{n\times n}$.

One such matrix operator,
$$
A_w=
\nabla^2 w-\frac { |\nabla w|^2 }{ 2w}I,
$$
where $I$ denotes the $n\times n$ identity matrix, has
arisen in conformal geometry
(see e.g. \cite{LL1}, \cite{V}  and the references therein).
In particular some comparision principles
for this matrix operator have been studied in 
\cite{L1} and \cite{L2}.

In almost all
 the papers mentioned above, 
 instead of supposing that  $\widetilde A_u$
lies in some region in ${\cal S}^{n\times n}$ it was
required that   the eigenvalues of 
$
\widetilde A_u$ lie in some
region $\Gamma$ in $\Bbb R^n$ which is symmetric,
i.e., if $\lambda =(\lambda_1, \cdots, \lambda_n)$
lies in $\Gamma$ then so does any permutation of the $\lambda_i$.
In this section, however, we will follow the more general 
approach.

We consider an open set ${\cal G}$ in ${\cal S}^{n\times n}$,
satisfying 
\begin{equation}
0\in \partial {\cal G},\ \ 
{\cal G}+{\cal P}\subset {\cal G},\ \
t{\cal G}\subset {\cal G}\ \forall\ t>0.
\label{2.6}
\end{equation}
Here ${\cal P}$ is the set of nonnegative matrices.

Consider functions $u, v\in C^2(B_1)\cap C^0(\overline B_1)$,
$B_1$ is the unit ball, satisfying (\ref{2.1}) and (\ref{2.2}).  In place of 
(\ref{2.3}) we require that
\begin{equation}
\widetilde A_u\in {\cal S}^{n\times n}\setminus {\cal G}\qquad\forall\ x\in B_1,
\label{2.7}
\end{equation}
\begin{equation}
\widetilde A_v\in \overline {\cal G}\qquad
\mbox{in}\ B_1.
\label{2.8}
\end{equation}
The generalized Hopf Lemma would be to conclude,
under possibly further conditions, that (\ref{2.4}) 
holds.

In the first result below we will assume 
\begin{equation}
u>0\qquad\mbox{in}\ B_1
\label{2.9}
\end{equation}
and
\begin{equation}
L(x, \beta s, \beta p)\le \beta L(x,s,p),
\qquad\forall\ \beta\ge 1,\ x,p\in \Bbb R^n, 
\ s>0.
\label{2.10}
\end{equation}

In order
to conclude that (\ref{2.4}) holds we impose, however, an additional condition
on ${\cal G}$.  It depends on the inner normal $\nu $ to $\partial B_1$ at $y$;
in our case $\nu=-y$; it is analogous to (\ref{2.5}).

\noindent\underline{Condition ${\cal G}_\nu$}: \
Let $\nu$ be a unit vector in $\Bbb R^n$. ${\cal G}$ is said to 
satisfy condition  ${\cal G}_\nu$ if there exists some open half 
cone $C_\delta(\nu)$,
$$
C_\delta(\nu)=\{ t(\{\nu_i\nu_j\}+A)\ |\
0<t, A\in {\cal S}^{n\times n},\
\|A\|<\delta\},\quad \delta>0
$$
such that
\begin{equation}
{\cal G}+C_\delta(\nu)\subset {\cal G}.
\label{2.11}
\end{equation}

Condition ${\cal G}_\nu$ cannot be simply dropped.  
In fact here is an example showing
that if ${\cal G}$ does not contain
$\{\nu_i\nu_j\}$ then (\ref{2.4})
fails in general.

\bigskip

\noindent\underline {Example 1}\ Assume that ${\cal G}$
satisfies  (\ref{2.6}) and  does not contain
$\{\nu_i\nu_j\}$, with $\nu=e_1=(1, 0, \cdots)$. Suppose 
$$
L(x, s, p)\le 0, \quad L(x, 1, 0)=0.
$$
Take
$$
u(x)=1+(x_1+1)^2, \quad v(x)\equiv 1.
$$
So $u>v$ in $\overline B_1\setminus\{-e_1\}$.  Here
$y=-e_1$
and $\nu=e_1$.

Then
$$
\widetilde A_v=0,
$$
$$
\widetilde A_u\le
 diag(2, 0, \cdots, 0), \mbox{therefore does not belong to}\ {\cal G},
$$
and (\ref{2.4}) does not hold.

\bigskip

\noindent\underline {Example 2}\
Assume that ${\cal G}$
satisfies  (\ref{2.6}) and $\{\nu_i\nu_j\}$, with $\nu=e_1$,
does not belong to $\overline {\cal G}$.
Assume
$$
L(x, 1, 0)\ge 0.
$$
Take, for large $k>1$,
$$
u(x)=1+k(x_1+1)^2, \quad v(x)\equiv 1.
$$
So $u>v$ in $\overline B_1\setminus\{-e_1\}$.  Here
$y=-e_1$
and $\nu=e_1$.

Then
$$
\widetilde A_v\ge 0,
$$
$$
\widetilde A_u=k[(2, 0, \cdots, 0)+O(\frac 1k)],\qquad\mbox{in}\
B_1\cap \{x\ |\ |x_1+1|k<1\}.
$$
So, for large $k$, 
$$
\widetilde A_u\ \mbox{does not belong to}\ \overline{\cal G},
$$
and (\ref{2.4}) does not hold.

\bigskip

Our first extension of the Hopf Lemma in this setup is

\begin{thm} Assume $u$ and $v$ satisfy (\ref{2.1}), (\ref{2.2}),
(\ref{2.9}), i.e. $u>0$ in
$B_1$, and (\ref{2.7}) and (\ref{2.8}), and that $L$ satisfies
(\ref{2.10}).  Assume furthermore
that ${\cal G}$ satisfies (\ref{2.6}) and condition ${\cal G}_\nu$, with
$\nu=-y$.  Then (\ref{2.4}) holds.
\label{theorem2.1}
\end{thm}

Here is an easy consequence of the theorem.

\noindent\underline{Corollary 1}.\ 
Let  $\Omega\subset \Bbb R^n$ be a bounded
domain with $C^2$ boundary,
 ${\cal G}$ satisfy (\ref{2.6}),
$$
diag(1, 0, \cdots, 0)\in  {\cal G},\quad
 {\cal G}\ \mbox{is convex},
$$
and
$$
O^t {\cal G} O\subset {\cal G}, \qquad \forall\ O\subset O(n),
$$
where $O(n)$ denotes the set of $n\times n$ orthoganal matrices,
and let $L$ satisfy (\ref{2.10}).
Suppose, for
 some $y\in \partial\Omega$,  that 
$u, v\in C^2(\Omega)\cap C^0(\Omega\cup \{y\})$,
 and satisfy
(\ref{2.1}) and (\ref{2.2}).
Then (\ref{2.4}) holds.

\medskip

\noindent\underline{Proof
of Theorem \ref{theorem2.1}:}\  Working with a smaller ball in
$B_1$ which contains $y$ on its boundary,
we may assume without loss of generality that $u>0$ in
$\overline B_1\setminus\{y\}$.
As usual, the proof makes use of
a comparison function:  Let
$$
h=e^{ -a|x|^2 }-e^{-a},
$$
where $a>1$ is a large constant to be chosen. 

For some  small constant, $\mu>0$, to be chosen, let
$$
\Sigma_\mu:= \{x\in B_1\ |\
-1<x\cdot (-y)< -1+\mu\}.
$$

 Differentiating $h$ we have
$$
h_i=-2 ax_i e^{ -a|x|^2}
$$
$$
h_{ij}=(4a^2 x_ix_j -2a \delta_{ij})  e^{ -a|x|^2}
=4a^2  e^{ -a|x|^2} (x_ix_j+O(\frac 1a) ).
$$

We will first fix the value of a small $\mu>0$
and  a large $a>1$,
and then fix the value of
 a small value $\epsilon>0$.

It follows that
$$
(v+\epsilon h)_{ij}=v_{ij}
+4a^2  e^{ -a|x|^2} \epsilon (x_ix_j+O(\frac 1a) ),
$$
and 
$$
L(x, v+\epsilon h, \nabla(v+\epsilon h)) 
-L(x, v, \nabla v)=\epsilon
O(h+a  e^{ -a|x|^2} )=\epsilon
O(a  e^{ -a|x|^2} ).
$$
Hence
$$
\widetilde A_{v+\epsilon h}=\widetilde A_v
+ 4a^2  e^{ -a|x|^2}  \epsilon (x_ix_j+O(\frac 1a) ).
$$

For small $\mu>0$ and large $a>1$, $x\in \Sigma_\mu$
 is close to $y$ and 
 the matrix $\{x_ix_j\}+O(\frac 1a)$ is close to
the matrix $\{\nu_i\nu_j\}$  and so lies in a cone
$C_\delta(\nu)$ as above.

Fixing the values of $\mu$ and $a$, 
then for all $x\in \Sigma_\mu$,
$$
 4a^2  e^{ -a|x|^2}  (x_ix_j+O(\frac 1a) )
$$
lies in the cone $C_\delta(\nu)$.  By condition ${\cal G}_\nu$, it follows
that 
$$
\widetilde A_{ v+\epsilon h}\in {\cal G}.
$$
Next, fix $0<\epsilon$ small so that
$$
u\ge v+\epsilon h\ \mbox{on}\ \partial \Sigma_\mu.
$$

\medskip

\noindent\underline{Claim}.\ $u\ge v+\epsilon h\ \mbox{in}\  \Sigma_\mu$.

\medskip

Indeed if the claim does not hold there exists
some constant $\beta>1$ and some 
$\bar x\in \Sigma_\mu$ such that
\begin{equation}
\beta u\ge v+\epsilon h\ \ \mbox{on}\
\Sigma_\mu,\qquad
\beta u(\bar x)=(v+\epsilon h)(\bar x).
\label{2.12}
\end{equation}
Here we have used the fact that $u>0$.

It follows, by (\ref{2.10}) that
$$
\beta \widetilde A_u(\bar x)\ge \widetilde A_{\beta u}(\bar x)
\ge  \widetilde A_{v+\epsilon h}(\bar x).
$$
But $ \widetilde A_{v+\epsilon h}(\bar x)\in {\cal G}$, hence
$$
\beta  \widetilde A_u(\bar x) \in {\cal G}
$$
and, by (\ref{2.6}), $
 \widetilde A_u(\bar x) \in {\cal G}$.  Contradiction.
The desired conclusion (\ref{2.4}) follows
from the claim.

\vskip 5pt
\hfill $\Box$
\vskip 5pt

Here is a slight variant of Theorem \ref{theorem2.1}.  In place of
(\ref{2.10}) we assume
\begin{equation}
L(x,s,p)\ \mbox{is nonincreasing in}\ s.
\label{2.10prime}
\end{equation}

\begin{thm}
Let $u, v$ and ${\cal G}$ satisfy the conditions of Theorem \ref{theorem2.1}
 except that we do not assume $u>0$, and we replace (\ref{2.10})
by (\ref{2.10prime}).  Then (\ref{2.4}) holds.
\label{theorem2.2}
\end{thm}

\noindent\underline{Proof}.\ The proof follows that of
Theorem \ref{theorem2.1} except that in the proof
of 
the claim, instead of choosing $\beta>1$ so that (\ref{2.12}) holds we choose 
a constant $\gamma>0$ so that
$$
u+\gamma \ge v+ \epsilon
h\ \mbox{in}\ \Sigma_\mu\ \ \mbox{and}\ \
u(\bar x)+\gamma =(v+\epsilon h)(\bar x).
$$
The rest of the argument is the same.

\vskip 5pt
\hfill $\Box$
\vskip 5pt

\bigskip

There is an open conjecture concerning a modified
kind of Hopf Lemma that arose in \cite{L-N2} and \cite{L-N3}, and to which
 we would like to call attention.

Consider $u\ge v$, positive functions
of $(t, y)$, $y\in \Bbb R^n$, in
$$
\Omega=\{(t, y)\ |\
0<t<1,\ |y|<1\}.
$$
and smooth in $\overline \Omega$.

Assume that
\begin{equation}
u_t>0\qquad \mbox{in}\ \Omega,
\label{2.13}
\end{equation}
\begin{equation}
u(0, y)\equiv 0\qquad\mbox{for}\ |y|<1,
\label{2.14}
\end{equation}
\begin{equation}
u_t(0, 0)=0,
\label{2.15}
\end{equation}
and the main condition:
$$
\mbox{whenever}\
u(t, y)=v(s, y), \ \mbox{for}\ 
0\le t\le s<1, \mbox{there} 
$$
\begin{equation}
\Delta u(t, y)\le \Delta v(s, y).
\end{equation}

\bigskip

\noindent\underline{Conjecture 1}.\ Then
\begin{equation}
u\equiv v.
\label{2.16}
\end{equation}

\bigskip

A weaker conjecture is

\bigskip

\noindent\underline{Conjecture 2}.\ Under the conditions above, (\ref{2.16}) holds
provided
$$
u(t, 0)\ \mbox{and}\ 
v(t, 0)\ \mbox{vanish of finite order
at}\ t=0.
$$

\bigskip

In \cite{L-N2} Conjecture 1 was proved in case
$u$ and $v$ are functions of $t$ alone;
\cite{L-N4} contains a simpler proof.
In the second paper, Conjecture 2 was proved in the 
very special case that, in addition,

\medskip

\noindent (a)\ \ the order of the first $t-$derivative of $u(t,0)$ which does not vanish
at the origin is $\le 3$.

\medskip

\noindent (b)\ \ $\nabla_y u_{tt}(0, 0)=0$.

\section{A remark on the Hopf Lemma for the heat operator}

We first recall the well known form; we describe it only for the 
classical heat operator but it holds, of course,
for much more general ones, see for example \cite{P-W}.

Consider a domain $G$ in $(x, t)$ space: here
$x\in \Bbb R^n$, $t\in \Bbb R$, lying in $t<0$, and
whose boundary includes an open domain $D$ on
$\{t=0\}$.  Assume the origin lies  in $\partial D$.
$\partial G\setminus D=:P\partial G$
is called the parabolic boundary of $G$.  For
convenience we suppose that
$(0, ...,  0, 1, 0)$ is the inner
normal to $\partial D$ at $(0,0)$, and
denote $x_n$ by $y$.  Sometimes we use
$(x,y,t)$ to denote a point, with
$x=(x_1, \cdots, x_{n-1})$.

Consider a function $u$ in 
$G\cup D\cup \{0,0\}$, $u>0$ in
$G\cup D$, $u(0, 0)=0$,
$$
u\in C^2(G)\cap C^0(G\cup D\cup \{0,0\}).
$$
Suppose 
\begin{equation}
(\partial _t -\Delta)u\ge 0\qquad\mbox{in}\ G.
\label{3.1}
\end{equation}
Here $\Delta$ is the Laplacian in the space variables, $(x,y)$.

For $P\partial G$ in $C^2$, the
parabolic
Hopf Lemma takes the form
\begin{thm}
Assume that the vector $(0, \cdots,0, 1, 0)$, the 
inner normal to $\partial D$ at $(0,0)$,  is not tangent to
$P\partial G$ at the origin.  Then
\begin{equation}
\liminf_{ s\to 0^+} \frac {  u(0, \cdots, 0, s, 0)  }s
>0.
\label{3.2}
\end{equation}
\label{theorem3.1}
\end{thm}

But something more may be true.  Suppose
that near the origin $D$ is
given by
$$
y>\alpha |x|^2, \ \ \alpha>0, \ \ t=0
$$
and assume that for some constant $a>0$, the domain
$$
\Omega=\{(x,y,t)\ |\
t<0, y>a|x|^2+c t\}
$$
near the origin, lies 
in $G$.

\begin{thm}
In this case we also have
\begin{equation}
\liminf_{t\to 0^-} \frac {  u(0, t)  }{ (-t) }>0.
\label{3.3}
\end{equation}
\label{theorem4.2}
\end{thm}

\noindent\underline{Proof}.\ We may suppose $u$ is continuous
in $\overline \Omega$ near the origin.  Indeed this may be achieved
by shrinking $\Omega$ slightly.  We use the comparison
function
$$
h=e^{  a(y-\alpha |x|^2 -ct) }.
$$
A calculation yields 
\begin{eqnarray*}
h_t-\Delta u&=&h
[-ac-a^2 +2a\alpha (n-1)
-4a^2 \alpha^2 |x|^2]\\
&<&0
\qquad\ \ \mbox{near the origin for}\ a\ \mbox{large}.
\end{eqnarray*}

For $0<R$ small consider the region
$$
\Omega_R=\Omega\cap \{|x|^2+y^2<R^2\}.
$$
On the parabolic boundary $P\partial \Omega_R$ we have
$$
u>0.
$$
Hence for $0<\epsilon$ small,
$$
u>\epsilon h\ \ \ \mbox{on}\ 
P\partial \Omega_R.
$$

From the parabolic maximum principle it follows that
\begin{equation}
u\ge \epsilon h\ \ \ \mbox{in}\
 \Omega_R.
\label{3.4}
\end{equation}

\vskip 5pt
\hfill $\Box$
\vskip 5pt

Since $h_t=-\alpha c$ at the origin, we infer from
(\ref{3.4}) that (\ref{3.3}) holds.

What happen if $(0, ..., 0, 1, 0)$ is tangent to 
$P\partial G$?
Consider the following simple example.  Suppose 
$G$ is given by 
$$
G=\{ (x,y, t)\ |\ t<0,
y-|x|^2 +\sqrt{-t}>0\}.
$$
(so $P\partial G$ is even analytic, namely,
$-t= (y-|x|^2)^2$)

\begin{thm} In this case, if $u$ satisfies the conditions
above, then
\begin{equation}
\liminf_{ t\to 0^-} \frac{  u(0, t) }{ \sqrt{-t} }>0.
\label{3.5}
\end{equation}
\label{theorem4.3}
\end{thm}
From this we see that $u(0, t)$ cannot be $C^1$
at the origin, or even $C^\alpha$ for $\alpha>1/2$.  In case of 
one space dimension there are known results on loss of regularity
when $(1,0)$ is tangent to $P\partial G$, see
\cite{K-N} and \cite{D}.  Our result
is for higher dimension and seems to exhibit
a phenomena not previously observed.  The loss of 
regularity is not due to the boundary values of $u$ near $0$, for there, $u$ could 
be $\equiv 0$, but still 
with $u$ positive in $G$.

\medskip

\noindent{\underline {Proof of Theorem \ref{theorem4.3}}}.\
For convenience we suppose $u$ is continuous in $\overline G$ near the 
origin
 (this may be achieved by replacing $G$ by 
$y-|x|^2+\sqrt{-\lambda t}>0$ with $0<1-\lambda$ small).
We makes use of the comparison function
$
h=y-|x|^2 +\sqrt{ -t}.
$
We have
$$
(\partial _t -\Delta)h=-\frac 1{  \sqrt{-t}}
+2(n-1)<0\ \ \mbox{for}\
|t|\ \mbox{small}.
$$ 
Arguing as above, we find that 
for some 
$0<\epsilon$ small,
$$
u\ge \epsilon h
$$
near the origin, and (\ref{3.5}) then follows.

\end{document}